\documentclass[10pt,reqno]{amsart}

\usepackage{a4wide}
\usepackage[active]{srcltx} 
\usepackage{amsmath,amssymb}
\usepackage[all]{xy}
\usepackage{mathbbol}
\usepackage{mathrsfs}
\usepackage{hyperref}

\newtheorem{proposition}{Proposition}[section]
\newtheorem{theorem}[proposition]{Theorem}
\newtheorem{lemma}[proposition]{Lemma}
\newtheorem{corollary}[proposition]{Corollary}

\theoremstyle{remark}
\newtheorem{remark}[proposition]{Remark}

\theoremstyle{definition}
\newtheorem{definition}[proposition]{Definition}

\newcommand{\cst}{\ifmmode\mathrm{C}^*\else{$\mathrm{C}^*$}\fi}
\newcommand{\st}{\;\vline\;}
\newcommand{\tens}{\otimes}
\newcommand{\atens}{\otimes_{\text{\rm\tiny{alg}}}}
\newcommand{\id}{\mathrm{id}}
\newcommand{\comp}{\!\circ\!}
\newcommand{\I}{\mathbb{1}}
\newcommand{\cpl}{\complement}
\newcommand{\ph}{\varphi}
\newcommand{\eps}{\varepsilon}
\newcommand{\CC}{\mathbb{C}}
\newcommand{\NN}{\mathbb{N}}
\newcommand{\is}[2]{\left\langle#1\,\vline\,#2\right\rangle}
\newcommand{\bra}[1]{\left\langle#1\,\vline\right.}
\newcommand{\ket}[1]{\left.\vline\,#1\right\rangle}
\newcommand{\GG}{\mathbb{G}}
\newcommand{\hGG}{\widehat{\mathbb{G}}}
\newcommand{\hDelta}{\widehat{\Delta}}
\newcommand{\GGmax}{\GG_{\mathrm{max}}}
\newcommand{\Deltamax}{\Delta_{\mathrm{max}}}
\newcommand{\GGmin}{\GG_{\mathrm{min}}}
\newcommand{\cc}{\mathrm{c}}
\newcommand{\op}{\mathrm{op}}
\newcommand{\ww}{\mathbb{w}}
\newcommand{\WW}{\mathbb{W}}
\newcommand{\VV}{\mathbb{V}}
\newcommand{\tp}{\xymatrix{*+<.7ex>[o][F-]{\scriptstyle\top}}}
\renewcommand{\Hat}[1]{\widehat{#1}}
\renewcommand{\Bar}[1]{\overline{#1}}
\DeclareMathOperator{\c0}{c_0}
\DeclareMathOperator{\Ad}{Ad}
\DeclareMathOperator{\M}{M}
\DeclareMathOperator{\HS}{HS}
\DeclareMathOperator{\B}{B}
\DeclareMathOperator{\C}{C}
\DeclareMathOperator{\cK}{\mathscr{K}}
\DeclareMathOperator{\cS}{\mathcal{S}}
\DeclareMathOperator{\Spec}{Spec}
\DeclareMathOperator{\ran}{ran}
\DeclareMathOperator{\Pol}{Pol}
\DeclareMathOperator{\Irr}{Irr}

\numberwithin{equation}{section}

\begin{document}
\title{Property $(\mathrm{T})$ and exotic quantum group norms}

\date{\today}

\author{David Kyed}
\address{Mathematisches Institut, Georg-August-Universit\"at G\"ottingen}
\email{kyed@uni-math.gwdg.de}
\urladdr{http://www.uni-math.gwdg.de/kyed/}

\author{Piotr M.~So{\l}tan}
\address{Institute of Mathematics, Polish Academy of Sciences\newline
\indent{and}\newline
\indent{}Department of Mathematical Methods in Physics, Faculty of Physics, University of Warsaw}
\email{piotr.soltan@fuw.edu.pl}
\urladdr{http://www.fuw.edu.pl/~psoltan/en/}

\thanks{P.S.~partially supported by Polish government grant no.~N201 1770 33, European Union grant PIRSES-GA-2008-230836 and Polish government matching grant no.~1261/7.PR UE/2009/7. D.K.~supported by The Danish Council for Independent Research $\mid$ Natural Sciences}

\dedicatory{Dedicated to Ryszard Nest on the occasion of his 60th birthday}

\begin{abstract}
Utilizing the notion of property $(\mathrm{T})$ we construct new examples of quantum group norms on the polynomial algebra of a compact quantum group, and provide criteria ensuring that these are not equal to neither the minimal nor the maximal norm. Along the way we generalize several classical operator algebraic characterizations of property $(\mathrm{T})$ to the quantum group setting which unify recent approaches to property $(\mathrm{T})$ for quantum groups with previous ones. The techniques developed furthermore provide tools to answer two open problems; firstly a question by B\'edos, Murphy and Tuset about automatic continuity of the comultiplication and secondly a problem left open by Woronowicz regarding the structure of elements whose coproduct is a finite sum of simple tensors.
\end{abstract}


\maketitle

\section{Introduction}

This paper is devoted to the theory of compact and discrete quantum groups. Both of these classes of quantum groups have been studied in detail by many authors and suffer from no shortage of exciting examples (\cite{su2,pw,podles,dqg,vawa}). It is known that a given compact quantum group $\GG$ can be described by more than one \cst-algebra (see e.g.~\cite{pseudogr,BMT01}); the most useful choices being the ``maximal'' and the ``minimal'' (also called \emph{reduced}) completions of the algebra $\Pol(\GG)$ of polynomial functions on $\GG$. It often happens that the canonical quotient map from the maximal completion $\C(\GGmax)$ to the minimal one $\C(\GGmin)$ is an isomorphism (in other words $\GG$ is \emph{co-amenable}, \cite{BMT01}). However, many interesting situations can arise when $\GG$ is described by a \cst-algebra sitting ``in between'' the maximal and minimal one (cf.~\cite{Dz}), but unfortunately there are not many examples of such compact quantum groups (apart from obvious direct product constructions).

Unlike compact quantum groups, the discrete quantum groups (i.e.~the duals of compact quantum groups) are all co-amenable --- there is just one \cst-algebra for each discrete quantum group. However, within this class of quantum groups one can find very interesting examples. In particular, there are discrete quantum groups with property $(\mathrm{T})$ which we will study in this paper. We will use property $(\mathrm{T})$ to construct special \cst-norms on the algebras of polynomials on (compact) dual quantum groups of property $(\mathrm{T})$ discrete quantum groups. The completions of these polynomial algebras will be ``exotic'' in the sense that they will sit in between the maximal and minimal completions. The canonical bijection between corepresentations of a discrete quantum group $\hGG$ and $*$-representations of the \cst-algebra $\C(\GGmax)$ will play a very important part in our investigation.

Let us briefly discuss the content of the paper. In Subsections \ref{nota} and \ref{prelim} we introduce the notation and list certain preliminary results from the theory of compact and discrete quantum groups. Section \ref{corep} provides necessary definitions and facts from the theory of corepresentations of quantum groups; we describe
the standard operations of forming tensor products and contragredient corepresentations, emphasizing the link with representations of the dual object. The regular corepresentation of a discrete quantum group is introduced in Subsection \ref{lambda} and in Subsection \ref{thm} we prove a quantum group version of a classical theorem from representation theory of locally compact groups. This theorem will be useful in the following sections.

In Section \ref{PT} property $(\mathrm{T})$ for discrete quantum groups is recalled and several known facts about quantum groups with property $(\mathrm{T})$ are listed. Then in Section \ref{Jacobson} the classical characterization of property $(\mathrm{T})$ in terms of isolated points in the space of irreducible representations is extended to the quantum group setting. This result provides a direct link between property $(\mathrm{T})$ of Fima (\cite{Fim08}) and earlier definitions in \cite{BCT05,Pet92}. As a consequence, in Section \ref{bb} we are able to show that a discrete quantum group $\hGG$ has property $(\mathrm{T})$ if and only if the \cst-algebra $\C(\GGmax)$ has property $(\mathrm{T})$ in the sense of Bekka (\cite{BekkaC}). Finally in Section \ref{minim} we extend the characterization of property $(\mathrm{T})$ by existence of a minimal projection in the full group \cst-algebra (\cite{Val84}) to the setting of discrete quantum groups.

The notion of a quantum group norm on the algebra of polynomials on a compact quantum group is defined in Section \ref{qgn}, where we also recall some basic facts about such norms. Then in Section \ref{adj} we give the first construction of a quantum group norm making the counit continuous. We call this procedure ``adjoining the neutral element to a compact quantum group.'' This construction provides examples of quantum group norms which differ from the reduced one as well as the maximal one in the absence of both amenability and property $(\mathrm{T})$. More complicated examples are collected in Section \ref{exo}, where, starting from a property $(\mathrm{T})$ discrete quantum group $\hGG$, we construct a certain quantum group norm $\|\cdot\|_\Pi$ on $\Pol(\GG)$. The completion of $\Pol(\GG)$ in this norm provides many examples of interesting (exotic) compact quantum groups. These examples lead us to a (negative) answer to a question of B\'edos, Murphy and Tuset whether any \cst-norm on $\Pol(\GG)$ arising from a representation weakly containing the regular one is necessarily a quantum group norm (cf.~\cite{BMT01} and Section \ref{exo}). Along the way we also give an example which shows that a very useful theorem of S.L.~Woronowicz about compact quantum groups with faithful Haar measure (\cite[Theorem 2.6(2)]{cqg}) cannot be generalized to all compact quantum groups. Several of our examples are co-commutative and we use some well known results from harmonic analysis (for which we refer e.g.~to \cite{BdlHV08} and geometric group theory (\cite{dlH}) to analyze them.

The paper uses the standard language of quantum group theory on operator algebra level (\cite{unbo,cqg,kv0}). In particular for a \cst-algebra $A$ the symbol $\M(A)$ will denote the multiplier algebra of $A$. All Hilbert spaces we will consider will be separable and the inner products will be linear in the second variable. Similarly all \cst-algebras, except multiplier algebras, will be assumed to be separable and the tensor product of \cst-algebras will always be the spatial one. The term ``representation'' will always mean a $*$-representation.

\subsection{Notation}\label{nota}

We shall adopt the convention of e.g.~\cite{pw,Fim08,Kye10} and always look at discrete quantum groups as duals of compact quantum groups. Thus any discrete quantum group will be denoted by $\hGG$. The \cst-algebra of ``continuous functions on $\hGG$ vanishing at infinity'' will be denoted by $\c0(\hGG)$ and its comultiplication by $\hDelta$. Thus
\[
\hGG=\bigl(\c0(\hGG),\hDelta\bigr).
\]
The compact quantum group $\GG$ dual to $\hGG$ can be described via many different objects. The polynomial algebra of $\GG$, i.e.~the Hopf $*$-algebra spanned by matrix elements of finite dimensional corepresentations of $\GG$, will be denoted by $\Pol(\GG)$. The universal enveloping \cst-algebra of $\Pol(\GG)$, i.e.~its completion with respect to the maximal \cst-norm will be denoted by $\C(\GGmax)$. The Hilbert space obtained via GNS construction from the Haar measure $h$ of $\GG$ will be denoted by $L^2(\GG)$. The completion of $\Pol(\GG)$ in the norm coming from representing $\Pol(\GG)$ on $L^2(\GG)$ will be denoted by $\C(\GGmin)$. Each of the algebras $\Pol(\GG)$, $\C(\GGmax)$ and $\C(\GGmin)$ has its own comultiplication, but we will use the same symbol $\Delta$ for all of them.

The possible other completions of $\Pol(\GG)$ will be denoted by $\C(\GG)$ or $\C(\GG_\Box)$, where in the space reserved by ``$\Box$'' a symbol indicating the nature of the completion will be placed. For example, if we choose a faithful representation $\pi$ of $\Pol(\GG)$ on some Hilbert space then the resulting \cst-completion of $\Pol(\GG)$ will be written as $\C(\GG_\pi)$. In case the \cst-norm used to complete $\Pol(\GG)$ is a quantum group norm (see Section \ref{qgn}) the \cst-algebra $\C(\GG_\Box)$ will carry a comultiplication extending that of $\Pol(\GG)$ and we will continue to denote it by the symbol $\Delta$. The only exception to this will come in parts of Section \ref{exo}, where the distinction between comultiplications on different completions of $\Pol(\GG)$ will be necessary. The von Neumann algebra obtained as the bicommutant of $\C(\GGmin)$ in $\B\bigl(L^2(\GG)\bigr)$ will be denoted by $L^\infty(\GG)$.

The set of equivalence classes of irreducible corepresentations of $\GG$ will be denoted by $\Irr(\GG)$. For $\alpha\in\Irr(\GG)$ we will choose (and fix throughout the paper) a corepresentation $u^\alpha$ in the class $\alpha$. The dimension of $u^\alpha$ will be denoted by $n_\alpha$. Thus $u^\alpha$ is a unitary element of $M_{n_\alpha}(\CC)\tens\Pol(\GG)$. As $\Pol(\GG)$ naturally embeds into $\C(\GGmax)$ and $\C(\GGmin)$ (or any $\C(\GG_\Box)$ for that matter), we can regard $u^\alpha$ as element of $M_{n_\alpha}(\CC)\tens\C(\GGmax)$ or
$M_{n_\alpha}(\CC)\tens\C(\GGmin)$ etc.

Let us also recall that a discrete quantum group $\hGG$ is \emph{unimodular} if its left and right Haar measures coincide (cf.~\cite[Section 3]{pw}). This is equivalent to many different conditions (cf.~\cite[Theorem 2.5]{cqg}). The one we will use is that of $\GG$ being a compact quantum group of \emph{Kac type} which manifests itself in the fact that the antipode of $\GG$ is a $*$-anti-automorphism.

\subsection{Some preliminary results}\label{prelim}

Recall from \cite[Section 3]{pw}, \cite[Section 4]{cqg} that
\[
\c0(\hGG)=\bigoplus_{\alpha\in\Irr(\GG)}M_{n_\alpha}(\CC)
\]
and it naturally acts on $L^2(\GG)$ which is the GNS Hilbert space for the Haar measure $h$ of $\GG$. We have the decomposition
\[
L^2(\GG)=\bigoplus_{\alpha\in\Irr(\GG)}H^\alpha
\]
with $H^\alpha$ the subspace of $L^2(\GG)$ spanned by matrix elements of $u^\alpha$. The set
\begin{equation}\label{uab}
\bigl\{u^\alpha_{i,j}\st\alpha\in\Irr(\GG),\;i,j\in\{1,\dotsc,n_\alpha\}\bigr\}
\end{equation}
is not an orthonormal basis of $L^2(\GG)$ in general (cf.~the Peter-Weyl-Woronowicz relations in \cite[Section 7]{cqg}), but if necessary the representatives $(u^\alpha)$ of classes in $\Irr(\GG)$ can be chosen so that it is an orthogonal system (\cite[Proposition 2.1]{daws}). Also let us note that if $\GG$ is of Kac type, then the system $\bigl\{\sqrt{n_\alpha}u^\alpha_{i,j}\st\alpha\in\Irr(\GG),\;i,j=1,\dotsc,n_\alpha\bigr\}$ is an orthonormal basis of $L^2(\GG)$.

In \cite[Section 2]{pw} the \emph{universal bicharacter} describing the duality between $\GG$ and $\hGG$ was introduced. It is the element
\begin{equation}\label{ww}
\ww=\bigoplus_{\alpha\in\Irr(\GG)}u^\alpha
\end{equation}
of $\M\bigl(\c0(\hGG)\tens\C(\GGmax)\bigr)$. It is of great importance and we will use it throughout the paper.

The action of $\c0(\hGG)$ on $L^2(\GG)$ is described in detail e.g.~in \cite{cqg}. Interpreting \cite[Formula 5.3]{cqg} in accordance with our notation we obtain for $a\in\Pol(\GG)$ and $\xi\in\C(\GGmax)^*$ the formula
\[
\bigl((\id\tens\xi)\ww\bigr)a=(\id\tens\xi)\Delta(a),
\]
where we view $\Pol(\GG)$ as a dense subspace of $L^2(\GG)$. Let us fix $\alpha$ and $i,j\in\{1,\dotsc,n_\alpha\}$ and take for $\xi$ the functional satisfying
\[
\xi(u_{k,l}^\beta)=\delta_{\alpha,\beta}\delta_{i,k}\delta_{j,l}
\]
(for existence of such a $\xi$ cf.~\cite[Section 1]{podles}) and put $a=u_{r,l}^\alpha$. Then we have
\[
(\id\tens\xi)\ww=e_{i,j}^\alpha\in{M_{n_\alpha}(\CC)}\subset
\bigoplus_{\alpha\in\Irr(\GG)}M_{n_\alpha}(\CC)=\c0(\Hat{\GG})
\]
and
\[
(\id\tens\xi)\Delta(u^\alpha_{r,l})=(\id\tens\xi)\sum_{k=1}^{n_\alpha}u^\alpha_{r,k}\tens{u^\alpha_{k,l}}=
\sum_{k=1}^{n_\alpha}\xi(u^\alpha_{k,l})u^\alpha_{r,k}
=\sum_{k=1}^{n_\alpha}\delta_{i,k}\delta_{j,l}u^\alpha_{r,k}=\delta_{j,l}u_{r,i}^\alpha.
\]
Thus $e^\alpha_{i,j}$ acts on a basic element $u^\alpha_{r,l}$ of $H^\alpha$ as
\begin{equation}\label{eua}
e^\alpha_{i,j}\colon{u^\alpha_{r,l}}\longmapsto{\delta_{j,l}u^\alpha_{r,i}}.
\end{equation}
The lesson from this is that if $H$ is a Hilbert space and $m\in{M_{n_\alpha}(\CC)}\tens\B(H)\subset\M\bigl(\c0(\Hat{\GG})\tens\cK(H)\bigr)$ is a matrix of operators
\[
m=\begin{bmatrix}
m_{1,1}&\dots&m_{1,n_\alpha}\\
\vdots&\ddots&\vdots\\
m_{n_\alpha,1}&\dots&m_{n_\alpha,n_\alpha}
\end{bmatrix}
\]
then for $r,l=1,\dotsc,n_\alpha$ and any $\xi\in{H}$ we have
\[
m(u^\alpha_{r,l}\tens\xi)=\sum_{i,j=1}^{n_\alpha}e^\alpha_{i,j}u^\alpha_{r,l}\tens{m_{i,j}}\xi=
\sum_{i,j=1}^{n_\alpha}\delta_{j,l}u^\alpha_{r,i}\tens{m_{i,j}}\xi
=\sum_{i=1}^{n_\alpha}u^\alpha_{r,i}\tens{m_{i,l}}\xi
\]
In particular, if $m(\eta\tens\xi)=\eta\tens\xi$ for any $\eta\in{H^\alpha}$ then taking $\eta=u^\alpha_{r,l}$ yields
\[
u^\alpha_{r,l}\tens\xi=\sum_{i=1}^{n_\alpha}u^\alpha_{r,i}\tens{m_{i,l}}\xi
\]
so that
\[
m_{i,l}\xi=\delta_{i,l}\xi.
\]

\section{Corepresentations of discrete quantum groups}\label{corep}

In this section we collect the standard facts about corepresentations of discrete quantum groups. Most of what we write here applies to all locally compact quantum groups and possibly more general quantum groups (cf.~\cite{mu2}), but in what follows we will stick with discrete quantum groups.

Let $\hGG$ be a discrete quantum group. A \emph{unitary corepresentation} of $\hGG$ on a Hilbert space $H_U$ is a unitary $U\in\M\bigl(\c0(\hGG)\tens\cK(H)\bigr)$ such that
\[
(\hDelta\tens\id)U=U_{23}U_{13}.
\]
(One usually expects the right hand side of the above equation to read $U_{13}U_{23}$, but this is really not so much different because $U^*$ satisfies such an equation.) Let $\ww$ be the universal bicharacter describing the duality between $\hGG$ and $\GG$ (defined by \eqref{ww}). Then it is known (\cite[Section 5.1]{mu2}) that any corepresentation $U\in\M\bigl(\c0(\hGG)\tens\cK(H_U)\bigr)$ is of the form
\[
U=(\id\tens\pi_U)\ww,
\]
where $\pi_U$ is a (uniquely determined) representation of $\C(\GGmax)$ on the Hilbert space $H_U$. We will not consider non-unitary corepresentations.

Let $U$ be a corepresentation of $\hGG$. Then $U=(\id\tens\pi)\ww$ for some representation $\pi$ of $\C(\GGmax)$. Now for any $\alpha\in\Irr(\GG)$ we define
\[
U^\alpha=(\id\tens\pi)u^\alpha.
\]
This is sometimes called the \emph{$\alpha$-component} of $U$, but note that $U^\alpha$ it is nothing like a sub-corepresentation.

Let us now describe some operations on corepresentations.

\subsection{Tensor product}

Take two corepresentations
\[
\begin{split}
U&=(\id\tens\pi_U)\ww\in\M\bigl(\c0(\Hat{\GG})\tens\cK(H_U)\bigr),\\
V&=(\id\tens\pi_V)\ww\in\M\bigl(\c0(\Hat{\GG})\tens\cK(H_V)\bigr)
\end{split}
\]
of $\hGG$. The \emph{tensor product} $U\tp{V}$ of $U$ and $V$ is defined as $U_{12}V_{13}\in\M\bigl(\c0(\Hat{\GG})\tens\cK(H_U\tens{H_V})\bigr)$. Another way to view the tensor product is
\[
U\tp{V}=\bigl(\id\tens[(\pi_U\tens\pi_V)\comp\Delta]\bigr)\ww.
\]
Indeed, $(\id\tens\Delta)\ww=\ww_{12}\ww_{13}$.

\subsection{Contragredient representation}

If $H$ is a Hilbert space and $\Bar{H}$ the complex conjugate Hilbert space then we have the anti-isomorphism
\[
\top\colon\B(H)\ni{m}\longmapsto\top(m)=m^\top\in\B(\Bar{H})
\]
given by
\[
m^\top\Bar{x}=\Bar{m^*x}.
\]

Let $V=(\id\tens\pi_V)\ww\in\M\bigl(\c0(\hGG)\tens\cK(H_V)\bigr)$ be a corepresentation. The \emph{contragredient} representation $V^\cc$ of $V$ is defined as $V^{\Hat{R}\tens\top}=(\Hat{R}\tens\top)V\in\M\bigl(\c0(\Hat{\GG})\tens\cK(\Bar{H_V})\bigr)$ (cf.~\cite[Section 3]{mu2}), where $\Hat{R}$ is the unitary antipode (\cite[Theorem 1.5(4)]{mu}). Again there is another way to view $V^\cc$:
\[
V^\cc=(\id\tens\pi_V^\cc)\ww,
\]
where
\[
\pi_V^\cc=\top\comp\pi_V\comp{R}
\]
and $R$ is the unitary antipode of $\GG$. This can be seen from
\[
\bigl(\id\tens[\top\comp\pi_V\comp{R}]\bigr)\ww=
(\Hat{R}\tens\top)(\id\tens\pi_V)(\Hat{R}\tens{R})\ww=(\Hat{R}\tens\top)(\id\tens\pi_V)\ww=V^\cc
\]
because $(\Hat{R}\tens{R})\ww=\ww$ (\cite[Formula 5.34]{mu2}).

\subsection{Containment, weak containment, equivalence, etc.}

Since there is a one to one correspondence between corepresentations of $\hGG$ and representations of the \cst-algebra $\C(\GGmax)$ we can define the notions of containment, weak containment, equivalence and weak equivalence of corepresentations by the corresponding notions from representation theory of \cst-algebras (see e.g.~\cite{Dix77} or Section \ref{Jacobson}). We will write $U\leq{V}$ if $U$ is contained in (i.e.~is a sub-corepresentation of) $V$ in the sense that $\pi_U$ is a subrepresentation of $\pi_V$. Similarly we will write $U\preccurlyeq{V}$ if $\pi_U\preccurlyeq\pi_V$ (weak containment). Two corepresentations $U$ and $V$ are \emph{equivalent} if $\pi_U$ and $\pi_V$ are unitarily equivalent, while $U$ and $V$ are \emph{weakly equivalent} if $U\preccurlyeq{V}$ and $V\preccurlyeq{U}$. We have the following simple lemma:

\begin{lemma}\label{little}
Let $U,U_1,V$ and $V_1$ be corepresentations of a discrete quantum group $\hGG$. Then
\begin{enumerate}
\item if $U\leq{U_1}$ and $V\leq{V_1}$ then $U\tp{V}\leq{U_1}\tp{V_1}$,
\item if $U\leq{V}$ then $U^\cc\leq{V^\cc}$.
\end{enumerate}
\end{lemma}

\begin{remark}\label{bG}
Let $U$ be a finite dimensional corepresentation of a discrete quantum group $\hGG$, i.e.~$U\in\M\bigl(\c0(\hGG)\tens\cK(H_U)\bigr)$ and $\dim{H_U}=n<\infty$. Then, upon choosing an orthonormal basis in $H_U$, we can identify $U$ with an $n\times{n}$ unitary matrix of elements of $\M\bigl(\c0(\hGG)\bigr)$ which satisfy
\[
\hDelta(U_{i,j})=\sum_{k=1}^nU_{k,j}\tens{U_{i,k}}
\]
for $i,j=1,\dotsc,n$. If we put $u_{i,j}=U_{j,i}^*$ then
\[
\hDelta(u_{i,j})=\sum_{k=1}^nu_{i,k}\tens{u_{k,j}}.
\]
Let $B$ be the \cst-subalgebra of $\M\bigl(\c0(\hGG)\bigr)$ generated by $\{u_{i,j}\}_{i,j=1,\dotsc,n}$. Then $B$ is unital (because $U$ is unitary) and $\hDelta$ restricts to a comultiplication $B\to{B}\tens{B}$. Then $\bigl(B,\bigl.\hDelta\bigr|_B\bigr)$ is a compact quantum matrix group (to see that condition 3.~of that definition is satisfied, consider the restriction of the antipode of $\hGG$ to the $*$-algebra generated by matrix elements of $U^*$, cf.~\cite[Theorem 1.6(4)]{mu}). Furthermore, $U$ is a unitary corepresentation of the opposite quantum group (\cite[Section 4]{kvvn}). Using the results of \cite[Section 3]{pseudogr}, \cite[Section 4]{kvvn} and \cite[Subsection 4.6]{qbohr} one can show that $U\tp{U^\cc}$ contains the trivial representation. Note, however, that the definition of contragredient corepresentation in \cite{pseudogr} is different from the one we have adopted and one is forced to use modular properties of the Haar measure of $\bigl(B,\bigl.\hDelta\!
 \bigr|_B\bigr)$.
\end{remark}

\subsection{The regular corepresentation}\label{lambda}

The regular corepresentation of a discrete quantum group $\hGG$ is $\WW=(\id\tens\lambda)\ww$, where $\lambda$ is the quotient map $\C(\GGmax)\to\C(\GGmin)\subset\B\bigl(L^2(\GG)\bigr)$.

\begin{proposition}\label{WW}
The regular corepresentation is equivalent to its contragredient $\WW^\cc$.
\end{proposition}

\begin{proof}
Let us first define a unitary map $Z\colon{L^2(\GG)}\to\Bar{L^2(\GG)}$. We put
\[
Zu^\alpha_{k,l}=\Bar{R({u^\alpha_{k,l}}^*)},
\]
where $R$ is the unitary antipode of $\GG$. The unitarity of $Z$ follows from the calculation:
\[
\begin{split}
\is{Zu^\alpha_{k,l}}{Zu^\beta_{i,j}}&=\is{R({u^\beta_{i,j}}^*)}{R({u^\alpha_{k,l}}^*)}\\
&=h\bigl(R({u^\beta_{i,j}}^*)^*R({u^\alpha_{k,l}}^*)\bigr)\\
&=h\bigl(R({u^\alpha_{k,l}}^*u^\beta_{i,j})\bigr)\\
&=h\bigl({u^\alpha_{k,l}}^*u^\beta_{i,j}\bigr)=\is{u^\alpha_{k,l}}{u^\beta_{i,j}}.
\end{split}
\]
Let us examine the operator $Z\lambda(a)Z^*$ for $a\in\Pol(\GG)$. On a vector $\Bar{R({u^\alpha_{k,l}}^*)}\in\Bar{L^2(\GG)}$ we have:
\[
\begin{split}
Z\lambda(a)Z^*\Bar{R({u^\alpha_{k,l}}^*)}&=Z\bigl(\lambda(a)u^\alpha_{k,l}\bigr)\\
&=Z(a\cdot{u^\alpha_{k,l}})\\
&=\Bar{R\bigl((au^\alpha_{k,l})^*\bigr)}\\
&=\Bar{R\bigl({u^\alpha_{k,l}}^*a^*\bigr)}\\
&=\Bar{R(a^*)\cdot{R}\bigl({u^\alpha_{k,l}}^*\bigr)}\\
&=\Bar{\lambda\bigl(R(a)\bigr)^*R\bigl({u^\alpha_{k,l}}^*\bigr)}\\
&=\lambda\bigl(R(a)\bigr)^\top\Bar{R\bigl({u^\alpha_{k,l}}^*\bigr)}.
\end{split}
\]
Thus $Z$ establishes unitary equivalence between $\lambda$ and $\top\comp\lambda\comp{R}$, which is the same as unitary equivalence between $\WW$ and $\WW^\cc$.
\end{proof}

\begin{remark}\label{regRep}
It is a well known fact that the tensor product $\WW\tp\WW$ is weakly contained in $\WW$ (cf.~\cite[Corollary 20]{mu2}). In view of Proposition \ref{WW}, we see that $\WW\tp\WW^\cc\preccurlyeq\WW$.
\end{remark}

\subsection{A theorem about corepresentations}\label{thm}

We end this section with a quantum group generalization of \cite[Proposition A.1.12]{BdlHV08} which gives a necessary and sufficient condition for a tensor product of two representations of a discrete group to have an invariant vector. Let $H$ and $K$ be Hilbert spaces and denote by $\HS(H,K)$ the space of Hilbert-Schmidt operators from $H$ to $K$. There is a canonical unitary mapping
\[
\Psi\colon{H}\tens{K}\longrightarrow\HS\bigl(\Bar{K},H\bigr)
\]
given by
\[
x\tens{y}\longmapsto\ket{x}\!\bra{\Bar{y}},
\]
where we use the Dirac notation: $\ket{x}\!\bra{\Bar{y}}$ is the operator
\[
\Bar{K}\ni\Bar{z}\longmapsto\is{\Bar{y}}{\Bar{z}}x\in{H}.
\]
This yields an isomorphism $\Ad_\Psi\colon\B(H\tens{K})\to\B\bigl(\HS\bigl(\Bar{K},H\bigr)\bigr)$
\[
\Ad_\Psi(x)=\Psi{x}\Psi^*.
\]

\begin{lemma}\label{maly}
For $S\in\B(H)$, $R\in\B(K)$ and $T\in\HS\bigl(\Bar{K},H\bigr)$ we have
\[
\bigl(\Ad_\Psi(S\tens{R})\bigr)(T)=S\comp{T}\comp{R^\top}.
\]
\end{lemma}

\begin{proof}
Calculate for $T=\Psi(x\tens{y})=\ket{x}\!\bra{\Bar{y}}$ and extend the result by linearity and continuity.
\end{proof}

\begin{theorem}\label{T}
Let $U$ and $V$ be corepresentations of a discrete quantum group $\hGG$. Then
\begin{enumerate}
\item\label{UVW1} if $W$ is a finite dimensional corepresentation of $\hGG$ such that $W\leq{U}$ and $W\leq{V^\cc}$ then $U\tp{V}$ contains the trivial corepresentation.
\item\label{UVW2} if $\hGG$ is unimodular and $U\tp{V}$ contains the trivial corepresentation then there exists a finite-dimensional corepresentation $W$ contained both in $U$ and in $V^\cc$.
\end{enumerate}
\end{theorem}

\begin{proof}
Ad \eqref{UVW1}. This follows directly from Lemma \ref{little} and Remark \ref{bG}.

Ad \eqref{UVW2}. As in the remarks preceding Lemma \ref{maly} we write $\Psi$ for the canonical unitary $H_U\tens{H_V}\to\HS\bigl(\Bar{H_V},H_U\bigr)$ and $\Ad_\Psi$ for the isomorphism $\B(H_U\tens{H_V})\to\B\bigl(\HS\bigl(\Bar{H_V},H_U\bigr)\bigr)$

Let us form the tensor product $U\tp{V}$. Let
\[
X=(\id\tens\Ad_{\Psi})(U\tp{V})\in\M\bigl(\c0(\hGG)\tens\cK\bigl(\HS\bigl(\Bar{H_V},H_U\bigr)\bigr)\bigr).
\]
Then $X$ is a corepresentation of $\hGG$ equivalent to $U\tp{V}$, so that $X$ contains the trivial corepresentation. This means that $X$ has a non-zero invariant vector.

Since
\[
X=(\id\tens\Ad_{\Psi})\bigl(\id\tens\bigl[(\pi_U\tens\pi_V)\comp\Delta\bigr]\bigr)\ww,
\]
the component $X^\alpha\in{M_{n_\alpha}(\CC)}\tens\B\bigl(\HS\bigl(\Bar{H_V},H_U\bigr)\bigr)$ is
\[
\begin{split}
X^\alpha&=(\id\tens\Ad_{\Psi})\bigl(\id\tens\bigl[(\pi_U\tens\pi_V)\comp\Delta\bigr]\bigr)u^\alpha\\
&=(\id\tens\Ad_{\Psi})\bigl(\id\tens\bigl[(\pi_U\tens\pi_V)\comp\Delta\bigr]\bigr)
\sum_{i,j=1}^{n_\alpha}e^\alpha_{i,j}\tens{u^\alpha_{i,j}}\\
&=(\id\tens\Ad_{\Psi})(\id\tens\pi_U\tens\pi_V)
\sum_{i,j,k=1}^{n_\alpha}e^\alpha_{i,j}\tens{u^\alpha_{i,k}}\tens{u^\alpha_{k,j}}\\
&=(\id\tens\Ad_{\Psi})\sum_{i,j}^{n_\alpha}e^\alpha_{i,j}\tens
\biggl(\sum_{k=1}^{n_\alpha}\pi_U(u^\alpha_{i,k})\tens\pi_V(u^\alpha_{k,j})\biggr)
\end{split}
\]
so that
\[
X^\alpha_{i,j}=\Ad_{\Psi}\biggl(\sum_{k=1}^{n_\alpha}\pi_U(u^\alpha_{i,k})\tens\pi_V(u^\alpha_{k,j})\biggr).
\]

By Lemma \ref{maly}, for $T\in\HS\bigr(\Bar{H_V},H_U\bigr)$ and $\eta\in{H^\alpha}$ we have
\[
X^\alpha(\eta\tens{T})=
\sum_{i,j}^{n_\alpha}e^\alpha_{i,j}\eta\tens
\biggl(\sum_{k=1}^{n_\alpha}\pi_U(u^\alpha_{i,k})\comp{T}\comp\pi_V(u^\alpha_{k,j})^\top\biggr)
\]

Now let $T$ be an invariant vector for $X$. In view of the discussion at the end of Subsection \ref{prelim}
\[
X^\alpha_{i,j}(T)=\delta_{i,j}T
\]
which reads
\[
\sum_{k=1}^{n_\alpha}\pi_U(u^\alpha_{i,k})\comp{T}\comp\pi_V(u^\alpha_{k,j})^\top
=\delta_{i,j}T.
\]
We have assumed that $\hGG$ is unimodular, i.e.~that $\GG$ is of Kac type. In particular, if $\kappa$ is the antipode of $\GG$ then $\kappa=R$ is a $*$-anti-automorphism and $\kappa^2=\id$. Moreover $\kappa(u^\alpha_{k,j})={u^\alpha_{j,k}}^*$. Therefore
\[
\sum_{k=1}^{n_\alpha}\pi_U(u^\alpha_{i,k})\comp{T}\comp\pi_V^\cc({u^\alpha_{j,k}}^*)
=\delta_{i,j}T
\]
Multiplying both sides of this equation by $\pi_V^\cc(u^\alpha_{j,p})$ and summing over $p$ we obtain
\[
\sum_{k=1}^{n_\alpha}\pi_U(u^\alpha_{i,k})\comp{T}\comp\pi_V^\cc\biggl(
\sum_{p=1}^{n_\alpha}{u^\alpha_{j,k}}^*u^\alpha_{j,p}
\biggr)
=\sum_{p=1}^{n_\alpha}\delta_{i,j}T\comp\pi_V^\cc(u^\alpha_{j,p})
\]
or equivalently
\begin{equation}\label{To}
\pi_U(u^\alpha_{i,p})\comp{T}=T\comp\pi_V^\cc(u^\alpha_{i,p})
\end{equation}
because $u^\alpha$ is unitary.

Since \eqref{To} is true for all $\alpha\in\Irr(\GG)$ and all $i,p\in\{1,\dotsc,n_\alpha\}$, we have
\[
U(\I\tens{T})=(\I\tens{T})V^\cc
\]
i.e.~$T$ intertwines $V^\cc$ and $U$.

It follows that $TT^*\in\cK(H_U)$ intertwines $U$ with itself. Note that $TT^*$ is a non-zero compact, positive operator. Therefore it has an eigenvalue $\lambda>0$ with finite multiplicity. Moreover the corresponding eigenprojection also intertwines $U$ with itself. This clearly leads to a finite dimensional sub-corepresentation $W$ of $U$. Similarly $T^*T$ is a self-intertwiner of $V^\cc$ and there is a sub-corepresentation $W'$ of $V^\cc$ corresponding to $\lambda$ (the non-zero parts of spectra of $TT^*$ and $T^*T$ coincide). Moreover, it is easy to see that the partial isometric part of the polar decomposition of $T^*$ establishes an equivalence between $W$ and $W'$.
\end{proof}

\begin{remark}
Let us remark that the first part of Theorem \ref{T} in the Kac case can be established in a simple calculation without resorting to the techniques described in Remark \ref{bG}. Indeed, using the notation of Theorem \ref{T} (and its proof), we first note that $U\tp{V}$ contains $W\tp{W^\cc}$. The corepresentation $W\tp{W^\cc}$ is equivalent to
\[
(\id\tens\widetilde{\pi})\ww\in\M\bigl(\c0(\hGG)\tens\cK(\HS(H_W,H_W))\bigr)
\]
where the representation $\widetilde{\pi}$ when restricted to $\Pol(\GGmax)$ is
\[
\widetilde{\pi}\colon\Pol(\GGmax)\ni{a}\longmapsto\Ad_\Psi\bigl((\pi_W\tens\pi_W)(\id\tens\kappa)\Delta(a)\bigr).
\]
In view of Lemma \ref{maly}, this means that for $T\in\HS(H_W,H_W)$ and $a\in\Pol(\GGmax)$ we have
\[
\bigl(\widetilde{\pi}(a)\bigr)(T)=\sum\pi_W(a_{(1)})\comp{T}\comp\pi_W\bigl(\kappa(a_{(2)})\bigr).
\]
Since $H_W$ is finite-dimensional, we can take $T=\I$ to obtain
\[
\bigl(\widetilde{\pi}(a)\bigr)(\I)=(\pi_W\tens\pi_W)\biggl(\sum(a_{(1)})\bigl(\kappa(a_{(2)})\bigr)\biggr)\I
=(\pi_W\tens\pi_W)\bigl(m(\id\tens\kappa)\Delta(a))\I=\eps(a)\I
\]
for all $a\in\Pol(\GGmax)$. It follows that the trivial corepresentation is contained in $W\tp{W^\cc}$.
\end{remark}

\section{\texorpdfstring{Property $(\mathrm{T})$ for discrete quantum groups}{Property (T) for discrete quantum groups}}\label{PT}

In a recent paper by P.~Fima \cite{Fim08}, Kazhdan's property $(\mathrm{T})$ is studied in the setting of discrete quantum groups. The definition is analogous to the classical definition for discrete groups and goes as follows.

\begin{definition}[{\cite{Fim08}}]\label{Tdef}
Let $\hGG$ be a discrete quantum group and let $V\in\M\bigl(\c0(\hGG)\tens\cK(H_V)\bigr)$ be a unitary corepresentation of $\hGG$ on the Hilbert space $H_V$. For $E\subset\Irr(\GG)$ a finite subset and $\delta>0$ a vector $\xi\in{H_V}$ is said to be \emph{$(E,\delta)$-invariant} with respect to $V$ if
\[
\bigl\|V^\alpha(\eta\tens\xi)-\eta\tens\xi\bigr\|<\delta\|\eta\|\|\xi\|
\]
for all $\alpha\in{E}$ and all $\eta\in{H^\alpha}$. The corepresentation $V$ \emph{has almost invariant vectors} if such a $\xi\in{H_V}$ exists for all finite subsets $E\subseteq\Irr(\GG)$ and all $\delta>0$, and the discrete quantum group $\hGG$ is said to have property $(\mathrm{T})$ if every corepresentation with almost invariant vectors has a non-zero invariant vector.
\end{definition}

\begin{remark}
It was shown in \cite{Kye10} how property $(\mathrm{T})$ for $\hGG$ can be interpreted using the correspondence between corepresentations of $\hGG$ and representations of $\C(\GGmax)$. More precisely, $\hGG$ has property $(\mathrm{T})$ if and only if the following holds: if $\pi\colon\C(\GGmax)\to\B(H)$ is a representation and there exists a sequence $(\xi_n)_{n\in\NN}$ of unit vectors $H$ such that $\lim\limits_{n\to\infty}\bigl\|\pi(a)\xi_n-\eps(a)\xi_n\bigr\|=0$ for all $a\in\C(\GGmax)$ then there exists a unit vector $\xi\in{H}$ with $\pi(a)\xi=\eps(a)\xi$ for all $a\in\C(\GGmax)$.
\end{remark}

\begin{remark}
Actually, the study of property $(\mathrm{T})$ for quantum groups began before the paper \cite{Fim08}. In \cite{Pet92} property $(\mathrm{T})$ was studied in the setting of Kac algebras and in \cite{BMT01} it was introduced for the class of algebraic quantum groups. However, we will use Fima's approach in the following, since it fits our purposes best. Using the results obtained in the present paper we will see later that the different approaches are equivalent in the case of discrete quantum  groups.
\end{remark}

In the following theorem we summarize some of the results obtained in \cite{Fim08}.
\begin{theorem}\label{ft}
If $\hGG$ is a discrete quantum group with property $(\mathrm{T})$ then the following holds:
\begin{enumerate}
\item $\hGG$ is finitely generated, i.e.~the compact dual is a matrix quantum group.
\item There exists a finite subset $E_0\subseteq\Irr(\GG)$ and a $\delta_0>0$ such that every corepresentation with $(E_0,\delta_0)$-invariant vectors has a non-zero invariant vector. Such a pair $(E_0,\delta_0)$ is called a \emph{Kazhdan pair} for $\hGG$.
\item\label{ftu} $\hGG$ is unimodular.
\end{enumerate}
\end{theorem}

Furthermore, Fima links property $(\mathrm{T})$ of $\hGG$ with property $(\mathrm{T})$ of $L^\infty(\GG)$ (in the sense of Connes and Jones \cite{cj}) in the case when $\hGG$ is i.c.c. Property $(\mathrm{T})$ for $\hGG$ can also be described by means of the ''positive definite functions'' on $\hGG$ as well as by a vanishing of cohomology result analogues to the classical Delorme-Guichardet theorem. We shall not elaborate further on these characterizations and refer the reader to \cite{Kye10} for details.

Property $(\mathrm{T})$ turns out to be essential in our search for exotic quantum group norms and in the following section we develop the results needed to construct these norms. The results obtained are of independent interests and parallel nicely classical results for discrete groups.

\section{\texorpdfstring{Property $(\mathrm{T})$ and the Jacobson topology}{Property (T) and the Jacobson topology}}\label{Jacobson}

Let again $\GG$ be a compact quantum group with \cst-algebra $\C(\GG)$ and Hopf $*$-algebra $\Pol(\GG)$. In this section we investigate the connection between property $(\mathrm{T})$ for $\hGG$ and the topology on the spectrum $\Spec\bigl(\C(\GGmax)\bigr)$ consisting of equivalence classes of irreducible representations of $\C(\GGmax)$. Recall from \cite{Dix77} that $\Spec\bigl(\C(\GGmax)\bigr)$ has a natural topology, called the Jacobson topology, which is intimately linked with the notion of weak containment. For the convenience of the reader we briefly recall this notion.

\begin{definition}[{\cite{Dix77}}]
Let $S$ be a set of representations of $\C(\GGmax)$ and $\pi$ some given representation. Then $\pi$ is said to be weakly contained in $S$, written $\pi\preccurlyeq{S}$, if every vector state associated with $\pi$ is a weak${}^*$ limit (i.e.~pointwise limit) of states which are linear combinations of vector functionals associated with the representations in $S$.
\end{definition}

\begin{proposition}[{\cite[Theorem 3.4.10]{Dix77}}]
If $S\subset\Spec\bigl(\C(\GGmax)\bigr)$ and $\pi\in\Spec\bigl(\C(\GGmax)\bigr)$ then the following are equivalent:
\begin{enumerate}
\item $\pi$ is in the closure of $S$ (with respect to the Jacobson topology) in $\Spec\bigl(\C(\GGmax)\bigr)$,
\item $\pi$ is weakly contained in $S$,
\item every vector state associated with $\pi$ is the weak${}^*$ limit of vector states associated with $S$.
\end{enumerate}
\end{proposition}

Recall that a functional $\ph\colon\C(\GGmax)\to\CC$ is said to be a vector functional associated with as set of representations $S$ if there exists $\rho\in{S}$ and $\xi\in{H_\rho}$ such that $\ph(a)=\is{\xi}{\rho(a)\xi}$ for all $a\in\C(\GGmax)$.

\begin{lemma}
Let $\pi\colon\C(\GGmax)\to\B(H)$ be a representation. Then $\pi$ has almost invariant vectors  if and only if the counit $\eps$ is weakly contained in $\pi$.
\end{lemma}

\begin{proof}
If $\pi$ has almost invariant vectors there exists a sequence $(\xi_n)_{n\in\NN}$ of unit vectors in $H$ such that
\[
\bigl\|\pi(a)\xi_n-\eps(a)\xi_n\bigr\|\xrightarrow[n\to\infty]{}0
\]
for all $a\in\C(\GG)$. Defining $\ph_n(a)=\is{\xi_n}{\pi(a)\xi_n}$ we have
\begin{equation}\label{one}
\bigl|\ph_n(a)-\eps(a)\bigr|^2=\bigl\|\pi(a)\xi_n-\eps(a)\xi_n\bigr\|^2-
\bigl(\ph_n(a^*a)-\ph_n(a^*)\ph_n(a)\bigr).
\end{equation}
and $\ph_n(a^*a)-\ph(a^*)\ph(a)\geq{0}$ by the Cauchy-Schwarz inequality. Thus
\[
\lim_{n\to\infty}\bigl|\ph_n(a)-\eps(a)\bigr|=0
\]
and we conclude that $\eps\preccurlyeq\pi$.

Conversely, if $\eps\preccurlyeq\pi$ we get a net $(\xi_\iota)$ of unit vectors in $H$ such that the net of vectors states $(\ph_\iota)$, $\ph_\iota(a)=\is{\xi_\iota}{\pi(a)\xi_\iota}$ converges pointwise to $\eps$ on $\C(\GGmax)$. But then for each $\iota$ we have 
\[
\bigl|\ph_\iota(a)-\eps(a)\bigr|^2=\bigl\|\pi(a)\xi_\iota-\eps(a)\xi_\iota\bigr\|^2-
\bigl(\ph_\iota(a^*a)-\ph_\iota(a^*)\ph_\iota(a)\bigr)
\]
as in \eqref{one} and hence $\lim\limits_{\iota}\bigl\|\pi(a)\xi_\iota-\eps(a)\xi_\iota\bigr\|=0$ for each $a\in\C(\GGmax)$. This shows that $\pi$ has almost invariant vectors.
\end{proof}

\begin{definition}
Let $\Omega$ be a set of positive functionals on $\C(\GGmax)$ and let $\ph$ be another positive functional. Then $\ph$ is said to be \emph{approximated on finite sets by elements in $\Omega$} if the following holds: for all finite $E\subset\Irr(\GG)$ and all $\delta>0$ there exists $\omega\in\Omega$ such that for all $\alpha\in{E}$ and all $i,j\in\{1,\dotsc,n_\alpha\}$
\[
\bigl|\ph(u^\alpha_{i,j})-\omega(u^\alpha_{i,j})\bigr|<\delta.
\]
\end{definition}

\begin{lemma}\label{onefive}
A positive functional $\ph$ on $\C(\GGmax)$ is approximated on finite sets by elements from $\Omega\subset\C(\GGmax)^*_+$ if and only if there exists a sequence $(\omega_n)_{n\in\NN}$ of elements of $\Omega$ such that $\omega_n(a)\xrightarrow[n\to\infty]{}\ph(a)$ for every $a\in\Pol(\GG)$. Moreover, in this case $\omega_n\xrightarrow[n\to\infty]{}\ph$ in the weak${}^*$ topology.
\end{lemma}

\begin{proof}
If $\ph$ is approximated by functionals from $\Omega$ on finite sets just pick an increasing sequence $(E_n)_{n\in\NN}$ of finite subsets of $\Irr(\GG)$ with $\Irr(\GG)$ as its union and choose $\omega_n\in\Omega$ such that
\[
\forall\;\alpha\in{E_n}\;\forall\;i,j\in\{1,\dotsc,n_\alpha\}\;\colon\;\bigl|\ph(u^\alpha_{i,j})-\omega_n(u^\alpha_{i,j})\bigr|<\tfrac{1}{n}.
\]
Since each matrix coefficient is contained in $E_n$ from a certain point on, we get the desired pointwise convergence on the set of matrix coefficients, and since these span $\Pol(\GG)$ linearly, the pointwise convergence holds on all of $\Pol(\GG)$. If, conversely, we have a sequence $(\omega_n)_{n\in\NN}$ of elements of $\Omega$ converging pointwise to $\ph$ on $\Pol(\GG)$ then clearly $\ph$ is approximated by functionals in $\Omega$ on finite subsets. That the convergence holds on all of $\C(\GGmax)$ is seen by a standard ``epsilon over three'' argument.
\end{proof}

\begin{remark}
Lemma \ref{onefive} shows that $\ph$ is approximated on finite sets by elements of $\Omega$ if and only if $\ph$ lies in the weak${}^*$ closure of $\Omega$.
\end{remark}

We will also use the following functional analytic version of a Lemma in \cite{BdlHV08}.

\begin{lemma}\label{oneseven}
If $\ph_1$ and $\ph$ are non-zero, positive linear functionals on $\C(\GGmax)$ and $\ph\geq\ph_1$ then the GNS representation $\pi_1$ associated with $\ph_1$ is contained in the GNS representation $\pi$ associated with $\ph$. If $\ph$ is already a vector functional associated to some representation $\rho$ then $\pi\leq\rho$.
\end{lemma}

\begin{proof}
Let $H$ and $H_1$ be the GNS Hilbert spaces associated to $\ph_1$ and $\ph$ respectively, with cyclic vectors $\Xi_1$ and $\Xi$. Since $\ph-\ph_1\geq{0}$ we get
\[
\bigl\|\pi_1(a)\Xi_1\bigr\|^2=\is{\Xi_1}{\pi_1(a^*a)\Xi_1}=\ph_1(a^*a)\leq\ph(a^*a)=\bigl\|\pi(a)\Xi\bigr\|^2.
\]
In particular
\[
\bigl(\pi(a)\Xi=0\bigr)\;\Longrightarrow\;\bigl(\pi_1(a)\Xi_1=0\bigr).
\]
Thus $T\colon{H}\to{H_1}$ defined by $\pi(a)\Xi\mapsto\pi_1(a)\Xi_1$ is well defined and bounded. Moreover, $T$ is trivially seen to intertwine $\pi$ and $\pi_1$. Put $K=(\ker{T})^\perp$ and note that $K$ is $\pi$-invariant. Since $T^*$ intertwines $\pi_1$ and $\pi$ the operator $T^*T$ is a self-intertwiner of $\pi$ and by functional calculus the
same is true for $|T|=(T^*T)^\frac{1}{2}$.

Consider now the polar decomposition $T=U|T|$, where $U$ is an isometry from $K$ onto $\Bar{\ran{T}}=H_1$. If we can prove that $U$ is also an intertwiner then $U^*$ provides us with an equivariant embedding of $H_1$ into $H$ proving that $\pi_1\leq\pi$. To see that $U$ intertwines we calculate for any $\xi\in{H}$
\[
\pi_1(a)U|T|\xi=\pi_1(a)T\xi=T\pi(a)\xi=U|T|\pi(a)\xi=U\pi(a)|T|\xi,
\]
which shows that $U$ restricted to $\Bar{\ran(|T|)}$ intertwines. But since $\Bar{\ran(|T|)}=\Bar{\ran(T^*)}=(\ker{T})^\perp=K$, we are done.

If furthermore $\ph(a)=\is{\eta}{\rho(a)\eta}$ for some representation $\rho$ on a Hilbert space $L$ then $V\colon{H}\to{L}$ given by $\pi(a)\Xi\mapsto\rho(a)\eta$ is a well defined isometry intertwining the GNS representation $\pi$ with $\rho$.
\end{proof}

\begin{remark}
Note that if (in the above proof) $\eta$ is cyclic for $\rho$ then $V$ is an equivalence between $\pi$ and $\rho$. This, for instance, is always the case if $\rho$ is irreducible.
\end{remark}

With the aid of the above lemmas we are now able to prove the following quantum group generalization of the classical characterization of property $(\mathrm{T})$ in terms of Fell's topology.

\begin{theorem}\label{isol}
A discrete quantum group $\hGG$ has property $(\mathrm{T})$ if and only if the trivial representation $\eps$ is an isolated point in $\Spec\bigl(\C(\GGmax)\bigr)$.
\end{theorem}

\begin{proof}
Assume first that $\hGG$ has property $(\mathrm{T})$. Since $\eps$ is finite dimensional and
irreducible $\{\eps\}$ is automatically closed in $\Spec\bigl(\C(\GGmax)\bigr)$, so we need to show that
$\{\eps\}$ is also open. Now the complement $\{\eps\}^\cpl$ is closed if and only if $\eps\not\in\Bar{\{\eps\}^\cpl}$ which happens if and only if $\eps$ is not weakly contained in $\{\eps\}^\cpl$. If this were the case then by \cite[Theorem 3.4.10]{Dix77} we find a net $(\pi_\iota)$ of elements of $\Spec\bigl(\C(\GGmax)\bigr)\setminus\{\eps\}$ and for each $\iota$ a unit vector $\xi_\iota$ in the representation space $H_\iota$ of $\pi_\iota$ such that the vector functionals $\ph_\iota\colon{a}\mapsto\is{\xi_\iota}{\pi_\iota(a)\xi_\iota}$ converge pointwise to $\eps$ on $\C(\GGmax)$. The functionals $\ph_\iota$ satisfy
\[
\bigl|\ph_\iota(a)-\eps(a)\bigr|^2=\bigl\|\pi_\iota(a)\xi_\iota-\eps(a)\xi_\iota\bigr\|^2-
\bigl(\ph_\iota(a^*a)-\ph_\iota(a^*)\ph_\iota(a)\bigr)
\]
(cf.~\eqref{one}) and hence $\bigl\|\pi_\iota(a)\xi_\iota-\eps(a)\xi_\iota\bigr\|\xrightarrow[\iota]{}0$ for all $a\in\C(\GGmax)$.

Define now
\[
\pi=\bigoplus_\iota\pi_\iota\colon\C(\GGmax)\longrightarrow\B\Bigl(\bigoplus_\iota{H_\iota}\Bigr).
\]
Then, by construction, $\pi$ has almost invariant vectors and by property $(\mathrm{T})$ it must have a non-zero invariant unit vector $\eta=(\eta_\iota)$. Then at least one $\eta_{\iota_0}$ is non-zero and hence invariant for $\pi_{\iota_0}$. Thus $\eps\leq\pi_{\iota_0}$ contradicting the choice of $\pi_{\iota_0}$ in $\Spec\bigl(\C(\GGmax)\bigr)\setminus\{\eps\}$.

Assume now that $\hGG$ does not have property $(\mathrm{T})$ and pick a representation $\pi\colon\C(\GGmax)\to\B(H)$ with almost invariant vectors but without non-zero invariant ones. We may therefore choose a sequence $(\xi_n)_{n\in\NN}$ of unit vectors in $H$ such that $\bigl\|\pi(a)\xi_n-\eps(a)\xi_n\bigr\|\xrightarrow[n\to\infty]{}0$ for every $a\in\Pol(\GG)$. Putting $\ph_n(a)=\is{\xi_n}{\pi(a)\xi_n}$ we obtain a sequence of positive functionals satisfying relation \eqref{one}. Hence $\ph_n(a)\xrightarrow[n\to\infty]{}\eps(a)$ for all $a\in\Pol(\GG)$ and Lemma \ref{onefive} assures that $\ph_n\to\eps$ in the weak${}^*$ topology. Our aim is to show that $\eps$ is not an isolated
point in the spectrum, i.e.~that $\eps$ is weakly contained in $\Spec\bigl(\C(\GGmax)\bigr)\setminus\{\eps\}$.

Hence, by Lemma \ref{onefive} we have to show that $\eps$ can be approximated on finite sets by elements from the set consisting of linear combinations of positive functionals associated with the representations in $\Spec\bigl(\C(\GGmax)\bigr)\setminus\{\eps\}$. Denote this set by $\Omega$.

Let $E$ be a finite subset of $\Irr(\GG)$ and let $\delta>0$ be given. Choose an $n_0\in{\NN}$ such that
\[
\bigl|\ph_{n_0}(u^\alpha_{i,j})-\eps(u^\alpha_{i,j})\bigr|<\tfrac{\delta}{2}
\]
for all $\alpha\in{E}$ and $i,j\in\{1,\dotsc,n_\alpha\}$. Recall also that $\ph_{n_0}(a)=\is{\xi_{n_0}}{\pi(a)\xi_{n_0}}$, with $\eps\nleq\pi$. Since the state space $\cS\bigl(\C(\GGmax)\bigr)$ is the weak${}^*$-closed convex hull of the the set of pure states of $\C(\GGmax)$, there exists a net $(\ph_\iota)$ of elements of $\cS\bigl(\C(\GGmax)\bigr)$ converging pointwise to $\ph_{n_0}$ and with the property that
\[
\ph_\iota=t_\iota\psi_\iota+(1-t_\iota)\eps,
\]
where $t_\iota\in[0,1]$ and $\psi_\iota$ is a linear combination of pure states different from $\eps$, i.e.~$\psi_\iota\in\Omega$. By compactness of $[0,1]$ and weak${}^*$ compactness of $\cS\bigl(\C(\GGmax)\bigr)$ we
may, upon passing to subnets, assume that $t_\iota\xrightarrow[\iota]{}t$ and $(\psi_{\iota})$ converges pointwise to a state $\psi$. Then
\[
\ph_{n_0}=t\psi+(1-t)\eps.
\]

If $t\neq{1}$ then Lemma \ref{oneseven} implies that $\eps$ is contained in the GNS representation associated to $\ph_{n_0}$ which, in turn, is contained in $\pi$ --- contradiction with the choice of $\pi$. Hence $t=1$ and thus $\ph_{n_0}$ is the pointwise limit of the net $(\psi_\iota)$. Hence there exists an index $\iota_0$ such that
\[
\bigl|\ph_{n_0}(u^\alpha_{i,j})-\psi_{\iota_0}(u^\alpha_{i,j})\bigr|<\tfrac{\delta}{2}
\]
for all $\alpha\in{E}$ and all $i,j\in\{1,\dotsc,n_\alpha\}$. Thus for all $\alpha\in{E}$ and $i,j\in\{1,\dotsc,n_\alpha\}$ we have $\bigl|\eps(u^\alpha_{i,j})-\psi_{\iota_0}(u^\alpha_{i,j})\bigr|<\delta$ and since $\psi_{\iota_0}$ is in the set $\Omega$, we have shown that $\eps$ is approximated by functionals in $\Omega$ as desired.
\end{proof}

\begin{remark}
Equipped with Theorem \ref{isol} one can easily prove that a discrete quantum group $\hGG$ has property $(\mathrm{T})$ in the sense of Definition \ref{Tdef} if and only if $\hGG$ has property $(\mathrm{T})$ as defined by E.~B\'edos, R.~Conti and L.~Tuset in \cite[Definition 7.15]{BCT05} and if and only if the associated Kac algebra (in von Neumann algebraic formulation) has property $(\mathrm{T})$ as defined by Petrescu and Joita in \cite[Definition 3.1]{Pet92} (cf.~\cite[Theorem 3.3]{Pet92}).
\end{remark}

\section{\texorpdfstring{Connection with property $(\mathrm{T})$ for \cst-algebras}{Connection with property (T) for C*-algebras}}\label{bb}

In the paper \cite{BekkaC} B.~Bekka introduced property $(\mathrm{T})$ for unital \cst-algebra endowed with a tracial state. His definition is a \cst-analogue of the corresponding definition for $\mathrm{II}_1$-factors due to
Connes and Jones (\cite{cj}) and goes as follows:

\begin{definition}
A unital \cst-algebra $A$ admitting a tracial state is said to have \emph{property $(\mathrm{T})$} if there exists a finite $F\subset{A}$ and a constant $c>0$ such that if a Hilbert $A$-bimodule $H$ has a unit vector $\xi$ such that
\[
\|a\xi-\xi{a}\|<c
\]
for all $a\in{F}$ then there exists a non-zero vector $\xi'\in{H}$ such that $a\xi=\xi{a}$ for all $a\in{A}$.
\end{definition}

\begin{theorem}\label{FiBe}
The discrete quantum group $\hGG$ has property $(\mathrm{T})$ if and only if $\C(\GGmax)$ has property $(\mathrm{T})$ in the sense of Bekka.
\end{theorem}

Note that the counit $\eps\colon\C(\GGmax)\to\CC$ is a tracial state so that Bekka's definition, which only covers \cst-algebras admitting tracial states, can be applied. The proof is greatly inspired by the the proof of \cite[Theorem 3]{Fim08}.

\begin{proof}[Proof of Theorem \ref{FiBe}]
Assume that $\hGG$ has property $(\mathrm{T})$ and let $(E,\delta)$ be a Kazhdan pair. We now prove that
\[
E'={\bigl\{u^\alpha_{i,j}\st\alpha\in{E},\;i,j\in\{1,\dotsc,n_\alpha\}\bigr\}}
\quad\text{and}\quad\delta'=\frac{\delta}{\max\{n_\alpha\sqrt{n_\alpha}\st\alpha\in{E}\}}
\]
constitute a Kazhdan pair for the \cst-algebra $\C(\GGmax)$. Assume therefore that $H$ is a Hilbert space which is also a $\C(\GGmax)$-bimodule and assume furthermore that $\xi$ is an $(E',\delta')$-central vector; i.e.~that
\[
\|u^\alpha_{i,j}\xi-\xi{u^\alpha_{i,j}}\|<\delta'
\]
for all $\alpha\in{E}$ and $i,j\in\{1,\dotsc,n_\alpha\}$. Denoting the left action $\C(\GGmax)\to\B(H)$ by $\pi$ and the right action $\C(\GGmax)^\op\to\B(H)$ by $\rho$ we obtain a new representation $\sigma\colon\C(\GGmax)\to\B(H)$ by setting $\sigma=m\comp\bigl([\rho\comp{R}]\tens\pi\bigr)\comp\Delta$ which corresponds to the corepresentation $V$ of $\hGG$ given by $V^\alpha=(\id\tens\rho)({u^\alpha}^*)(\id\tens\pi)(u^\alpha)$. For $\alpha\in{E}$ we now obtain, using \eqref{eua} and the fact that $\GG$ is Kac so that 
\begin{equation}\label{orth}
\bigl\{\sqrt{n_\alpha}u^\alpha_{i,j}\st{i,j=1,\dotsc,n_\alpha}\bigr\}
\end{equation}
is an orthonormal basis of $H^\alpha$, that
{\allowdisplaybreaks
\begin{align*}
\bigl\|V^\alpha(u^\alpha_{r,l}\tens\xi)-u^\alpha_{r,l}\tens\xi\bigr\|
&=\bigl\|\bigl((\id\tens\pi)(u^\alpha)\bigr)(u^\alpha_{r,l}\tens\xi)
-\bigl((\id\tens\rho)(u^\alpha)\bigr)(u^\alpha_{r,l}\tens\xi)\bigr\|\\
&=\biggl\|\sum_{i,j=1}^{n_\alpha}\bigl(e^\alpha_{i,j}\tens\pi(u^\alpha_{i,j})\bigr)(u^\alpha_{r,l}\tens\xi)
-\sum_{i,j=1}^{n_\alpha}\bigl(e^\alpha_{i,j}\tens\rho(u^\alpha_{i,j})\bigr)(u^\alpha_{r,l}\tens\xi)\biggr\|\\
&=\biggl\|\sum_{i=1}^{n_\alpha}u^\alpha_{r,i}\tens\pi(u^\alpha_{i,l})\xi
-\sum_{i=1}^{n_\alpha}u^\alpha_{r,i}\tens\rho(u^\alpha_{i,l})\xi\biggr\|\\
&=\biggl\|\sum_{i=1}^{n_\alpha}u^\alpha_{r,i}\tens(u^\alpha_{i,l}\xi-\xi{u^\alpha_{i,l}})\biggr\|\\
&=\sqrt{\sum_{i=1}^{n_\alpha}\|u_{r,i}^\alpha\|^2\| u^\alpha_{i,l}\xi-\xi{u^\alpha_{i,l}} \|^2}\\
&=\sqrt{\sum_{i=1}^{n_\alpha}\tfrac{1}{n_\alpha}\| u^\alpha_{i,l}\xi-\xi{u^\alpha_{i,l}} \|^2}
<\delta'.
\end{align*}} 

Now we take $\eta\in{H^\alpha}$ and expand it in the orthonormal basis \eqref{orth}:
\[
\eta=\sum_{r,i=1}^{n_\alpha}n_\alpha\is{u^\alpha_{r,i}}{\eta}u^\alpha_{r,i}.
\]
Then
\[
\begin{split}
\big\|V^\alpha(\eta\tens\xi)-\eta\tens\xi\bigr\|
&=\biggl\|\sum_{r,i=1}^{n_\alpha}n_\alpha\is{u^\alpha_{r,i}}{\eta}
\bigl(V^\alpha(u^\alpha_{r,l}\tens\xi)-u^\alpha_{r,l}\tens\xi\bigr)\biggr\|\\
&\leq\sum_{r,i=1}^{n_\alpha}n_\alpha\bigl|\is{u^\alpha_{r,i}}{\eta}\bigr|
\bigl\|V^\alpha(u^\alpha_{r,l}\tens\xi)-u^\alpha_{r,l}\tens\xi\bigr\|\\
&<\sum_{r,i=1}^{n_\alpha}n_\alpha\bigl|\is{u^\alpha_{r,i}}{\eta}\bigr|\delta'
=\delta'\sqrt{n_\alpha}\sum_{r,i=1}^{n_\alpha}\bigl|\is{\sqrt{n_\alpha}u^\alpha_{r,i}}{\eta}\bigr|
\end{split}
\]
By the Cauchy-Schwarz inequality
\[
\sum_{r,i=1}^{n_\alpha}\bigl|\is{\sqrt{n_\alpha}u^\alpha_{r,i}}{\eta}\bigr|
\leq\sqrt{\sum_{r,i=1}^{n_\alpha}\bigl|\is{\sqrt{n_\alpha}u^\alpha_{r,i}}{\eta}\bigr|^2}
\sqrt{\sum_{r,i=1}^{n_\alpha}1}=\|\eta\|n_\alpha,
\]
so that
\[
\big\|V^\alpha(\eta\tens\xi)-\eta\tens\xi\bigr\|
<\delta'\sqrt{n_\alpha}\sum_{r,i=1}^{n_\alpha}\bigl|\is{\sqrt{n_\alpha}u^\alpha_{r,i}}{\eta}\bigr|
\leq\delta'\sqrt{n_\alpha}\|\eta\|n_\alpha=\delta'\|\eta\|n_\alpha^{\frac{3}{2}}\leq\delta\|\eta\|
\]
for all $\eta\in{H^\alpha}$. Therefore, since $(E,\delta)$ is a Kazhdan pair for $\hGG$, there exists a $V$-invariant unit vector $\xi'\in{H}$. It is easily seen that $\xi'$ is a central vector and we conclude that $\C(\GGmax)$ has property $(\mathrm{T})$.

If, conversely, $\C(\GGmax)$ has property $(\mathrm{T})$ then from \cite[Proposition 3.2]{Bro06} it follows that every finite dimensional, irreducible representation is an isolated point in $\Spec\bigl(\C(\GGmax)\bigr)$. In particular $\eps$ is an isolated point and therefore $\hGG$ has property $(\mathrm{T})$ by Theorem \ref{isol}.
\end{proof}

\begin{remark}\label{fin}
Let us emphasize that Theorem \ref{FiBe} together with \cite[Proposition 3.2]{Bro06} shows that, as in the classical case, $\hGG$ has property $(\mathrm{T})$ if and only if all finite dimensional representations of $\C(\GGmax)$ are isolated in $\Spec\bigl(\C(\GGmax)\bigr)$.
\end{remark}

\begin{corollary}\label{notweakly}
Let $\hGG$ be an infinite discrete quantum group. i.e.~one with $\dim\c0(\hGG)=\infty$. Assume that $\hGG$ has property $(\mathrm{T})$. Then the regular corepresentation $\WW$ of $\hGG$ does not weakly contain any finite dimensional corepresentation.
\end{corollary}

\begin{proof}
If $U$ is a finite dimensional corepresentation of $\hGG$ and $U\preccurlyeq\WW$ then by Remark \ref{fin} we have $U\leq\WW$. By Lemma \ref{little} and Remark \ref{regRep} we have
\[
U\tp{U^\cc}\leq\WW\tp\WW^\cc\preccurlyeq\WW.
\]
But $U\tp{U^\cc}$ contains the trivial representation (Remark \ref{bG}), so $\WW$ must weakly contain the trivial representation which is impossible for property $(\mathrm{T})$ infinite discrete quantum group (cf.~\cite[Theorem 7.17]{BCT05}, \cite[Remark 4]{Fim08}).
\end{proof}

\section{\texorpdfstring{Minimal projections and property $(\mathrm{T})$}{Minimal projections and property (T)}}\label{minim}

We recall that a projection $p$ in a unital \cst-algebra $A$ is called \emph{minimal} if $pAp=\CC{p}$. We prove here the following quantum group version of the classical characterization of property $(\mathrm{T})$ in terms of minimal projections in the maximal group \cst-algebra (see \cite{aw,Val84}). The proof follows the lines of the corresponding proof in \cite{Val84}. As usual $\hGG$ denotes a discrete quantum group.

\begin{proposition}\label{miprojT}
The following are equivalent:
\begin{enumerate}
\item\label{21one} $\hGG$ has property $(\mathrm{T})$,
\item\label{21two} there exists a unique minimal projection in the center of $\C(\GGmax)$ with $\eps(p)=1$,
\item\label{21three} there exists a minimal projection $p\in\C(\GGmax)$ with $\eps(p)=1$.
\end{enumerate}
\end{proposition}

\begin{proof}
We first prove \eqref{21one}$\Rightarrow$\eqref{21two}. If $\hGG$ has property $(\mathrm{T})$ then $\eps$ is isolated in $\Spec\bigl(\C(\GGmax)\bigr)$. Hence the spectrum splits into a disjoint union of open subsets as
$\Spec\bigl(\C(\GGmax)\bigr)=\{\eps\}\cup\{\eps\}^\cpl$ and thus $\C(\GGmax)$ splits accordingly (as a \cst-algebra!)
into the direct sum of two closed, two-sided ideals $I$ and $J$ defined, implicitly, as
\[
\begin{split}
\{\eps\}&=\bigl\{\pi\in\Spec\bigl(\C(\GGmax)\bigr)\st\pi(I)\neq\{0\}\bigr\},\\
\{\eps\}^\cpl&=\bigl\{\pi\in\Spec\bigl(\C(\GGmax)\bigr)\st\pi(J)\neq\{0\}\bigr\}.
\end{split}
\]
Clearly we have $J=\ker{\eps}$ and thus $I$ is one-dimensional. The unit now splits as $\I=(e,f)$ and $p=(e,0)$ clearly does the job. If another minimal, central projection $p'$ with $\eps(p')=1$ existed then we would have $pp'=pp'p=\lambda{p}$ for some $\lambda\in\CC$ and since $\eps(p)=\eps(p')=1$ we have $\lambda=1$. Thus $p\leq{p'}$ and by minimality $p=p'$.

The implication \eqref{21two}$\Rightarrow$\eqref{21three} is obvious. Lastly we prove  \eqref{21three}$\Rightarrow$\eqref{21one}. Let therefore $p\in\C(\GGmax)$ be minimal with $\eps(p)=1$. Then by \cite[Lemma 1]{Val84} there exists a unique
$\pi\in\Spec\bigl(\C(\GGmax)\bigr)$ such that $\pi(p)\neq{0}$. Since $\eps$ clearly is such a representation we
have
\[
\{\eps\}=\bigl\{\rho\in\Spec\bigl(\C(\GGmax)\bigr)\st\rho(p)\neq{0}\bigr\}
=\bigl\{\rho\in\Spec\bigl(\C(\GGmax)\bigr)\st\rho\bigl(\C(\GGmax)p\C(\GGmax)\bigr)\neq\{0\}\bigr\}
\]
which by definition is open in $\Spec\bigl(\C(\GGmax)\bigr)$. Since $\{\eps\}$ is always closed, this proves that $\eps$ is an isolated point in $\Spec\bigl(\C(\GGmax)\bigr)$ and hence that $\hGG$ has property $(\mathrm{T})$.
\end{proof}

\section{Quantum group norms}\label{qgn}

In \cite[Section 3]{BMT01} the question of completing the polynomial algebra $\Pol(\GG)$ under different \cst-norms was addressed. In particular, the authors considered \cst-norms on $\Pol(\GG)$ for which the comultiplication $\Pol(\GG)\to\Pol(\GG)\atens\Pol(\GG)$ extends to a $*$-homomorphism of the completions. Such norms were called \emph{regular.} We feel that the term ``regular'' is already overused in the literature on quantum groups (let us mention e.g.~the \emph{regularity} condition for multiplicative unitaries of \cite{bs,bsv} or the \emph{regular} corepresentation of Subsection \ref{lambda}). Therefore we would like to propose the following terminology:

\begin{definition}
Let $\GG$ be a compact quantum group and let $\|\cdot\|_\sim$ be a \cst-norm on $\Pol(\GG)$. Let $\C(\GG_\sim)$ be the completion of $\Pol(\GG)$ in the norm $\|\cdot\|_\sim$. The \cst-norm $\|\cdot\|_\sim$ is called a \emph{quantum group norm} if the comultiplication $\Delta\colon\Pol(\GG)\to\Pol(\GG)\atens\Pol(\GG)$ extends to a $*$-homomorphism $\C(\GG_\sim)\to\C(\GG_\sim)\tens\C(\GG_\sim)$.
\end{definition}

B\'edos, Murphy and Tuset proved, among other things, that the norm coming from the representation of $\Pol(\GG)$ on $L^2(\GG)$ is the smallest quantum group norm on $\Pol(\GG)$ (cf.~Remark \ref{regRep} for an argument that it is a quantum group norm). Also the \emph{universal} or \emph{maximal} \cst-norm on $\Pol(\GG)$, i.e.~the supremum of all \cst-norms on $\Pol(\GG)$, was proved in \cite{BMT01} to be a quantum group norm.

In the next sections we will construct examples of quantum group norms with various interesting properties. In particular we will obtain examples of compact quantum groups $\GG$ sitting strictly ``between'' their minimal and maximal versions. We will provide such examples both admitting a continuous counit and without this property.

\section{Adjoining the neutral element to a compact quantum group}\label{adj}

Let $\GG$ be a compact quantum group. We may view $\C(\GG)$ as embedded into $\B(H)$ for some Hilbert space $H$, so that the inclusion $\Pol(\GG)\hookrightarrow\C(\GG)$ becomes a representation of the $*$-algebra $\Pol(\GG)$, say $\pi\colon\Pol(\GG)\to\B(H)$. Consider now the representation
\begin{equation}\label{witipi}
\widetilde{\pi}\colon\Pol(\GG)\ni{a}\longmapsto
\begin{bmatrix}
\pi(a)&0\\0&\eps(a)
\end{bmatrix}
\in\B(H\oplus\CC)
\end{equation}
and let $\|\cdot\|_{\widetilde{\pi}}$ be the norm defined by $\widetilde{\pi}$:
\begin{equation}\label{tipi}
\|a\|_{\widetilde{\pi}}=\bigl\|\widetilde{\pi}(a)\bigr\|
=\max\bigl\{\bigl\|\pi(a)\bigr\|,\bigl|\epsilon(a)\bigr|\bigr\}
\end{equation}
for all $a\in\Pol(\GG)$.

\begin{proposition}
The \cst-norm $\|\cdot\|_{\widetilde{\pi}}$ on $\Pol(\GG)$ is a quantum group norm.
\end{proposition}

\begin{proof}
Take $a\in\Pol(\GG)$. We have
\[
\begin{split}
(\widetilde{\pi}\tens&\widetilde{\pi})\Delta(a)
=\sum\widetilde{\pi}(a_{(1)})\tens\widetilde{\pi}(a_{(2)})\\
&=\sum
\begin{bmatrix}
\pi(a_{(1)})\tens\pi(a_{(2)})&0&0&0\\
0&\pi(a_{(1)})\eps(a_{(2)})&0&0\\
0&0&\eps(a_{(1)})\pi(a_{(2)})&0\\
0&0&0&\eps(a_{(1)})\eps(a_{(2)})
\end{bmatrix}\\
&=
\begin{bmatrix}
(\pi\tens\pi)\Delta(a)&0&0&0\\
0&\pi(a)&0&0\\
0&0&\pi(a)&0\\
0&0&0&\eps(a)
\end{bmatrix}.
\end{split}
\]
Therefore
\[
\bigl\|(\widetilde{\pi}\tens\widetilde{\pi})\Delta(a)\bigr\|
=\max\bigl\{\bigl\|(\pi\tens\pi)\Delta(a)\bigr\|,\bigl\|\pi(a)\bigr\|,\bigl|\eps(a)\bigr|\bigr\}.
\]
Since the norm defined by $\pi$ is a quantum group norm, we have $\bigl\|(\pi\tens\pi)\Delta(a)\bigr\|\leq\bigl\|\pi(a)\bigr\|$ for all $a\in\Pol(\GG)$. It follows that
\begin{equation}\label{ten}
\bigl\|(\widetilde{\pi}\tens\widetilde{\pi})\Delta(a)\bigr\|
=\max\bigl\{\bigl\|\pi(a)\bigr\|,\bigl|\eps(a)\bigr|\bigr\}=\bigl\|\widetilde{\pi}(a)\bigr\|.
\end{equation}
Let $\C(\GG_{\widetilde{\pi}})$ be the completion of $\Pol(\GG)$ with respect to the norm $\|\cdot\|_{\widetilde{\pi}}$. Then \eqref{ten} shows that $\Delta\colon\Pol(\GG)\to\Pol(\GG)\atens\Pol(\GG)$ extends to an isometry $\C(\GG_{\widetilde{\pi}})\to\C(\GG_{\widetilde{\pi}})\tens\C(\GG_{\widetilde{\pi}})$ (minimal -- spatial -- tensor product).
\end{proof}

\begin{definition}
Let $\GG$ be a compact quantum group. The compact quantum group obtained by the construction described in the above proposition will be denoted by $\widetilde{\GG}$ and called the quantum group $\GG$ \emph{with neutral element adjoined.} Thus, by definition $\C(\widetilde{\GG})=\C(\GG_{\widetilde{\pi}})$, where $\widetilde{\pi}$ is defined by \eqref{witipi}.
\end{definition}

\begin{proposition}\label{mini}
Assume that $\GG$ does not admit a continuous co-unit. Then there exists a central projection $p$ in $\C(\widetilde{\GG})$ such that
\[
\C(\widetilde{\GG})\cong\CC{p}\oplus\C(\GG).
\]
\end{proposition}

\begin{proof}
As before we write $\pi$ for the representation $\Pol(\GG)\hookrightarrow\C(\GG)\subset\B(H)$ for some Hilbert space $H$ and $\widetilde{\pi}$ for the direct sum of $\pi$ and $\eps$. Denote by $\|\cdot\|_\pi$ and $\|\cdot\|_{\widetilde{\pi}}$ the associated \cst-norms on $\Pol(\GG)$. Since $\eps$ is unbounded on $(\Pol(\GG),\|\cdot\|_\pi)$, for each $n\in\NN$ there exists $a_n\in\Pol(\GG)$ such that $\|a_n\|_\pi=1$ and $\bigl|\eps(a_n)\bigr|>n$. Let $b_n=\tfrac{1}{\eps(a_n)}a_n$. Clearly $\bigl\|\pi(b_n)\bigr\|=\|b_n\|_\pi\xrightarrow[n\to\infty]{}0$, while $\eps(b_n)=1$ for all $n$.

The completion $\C(\widetilde{\GG})$ of $\Pol(\GG)$ in $\|\cdot\|_{\widetilde{\pi}}$ is isomorphic to the closure of $\widetilde{\pi}\bigl(\Pol(\GG)\bigr)$ inside $\B(H\oplus\CC)$. Note that the sequence $\bigl(\widetilde{\pi}(b_n)\bigr)_{n\in\NN}$ converges in $\B(H\oplus\CC)$, since
\[
\widetilde{\pi}(b_n)=
\begin{bmatrix}
\pi(b_n)&0\\0&\eps(b_n)
\end{bmatrix}
\]
and clearly
\[
p=\lim_{n\to\infty}\widetilde{\pi}(b_n)=
\begin{bmatrix}
0&0\\0&1
\end{bmatrix}.
\]
It is now clear that $p$ commutes with all elements of the form
\[
\begin{bmatrix}
\pi(a)&0\\0&\eps(a)
\end{bmatrix}
\]
($a\in\Pol(\GG)$) and that we have the isomorphism $\C(\widetilde{\GG})\cong\CC{p}\oplus\C(\GG)$.
\end{proof}

\begin{remark}
\noindent
\begin{enumerate}
\item If $\GG$ is a compact quantum group with continuous counit (e.g.~$\GG$ might be co-amenable, cf.~\cite{BMT01}) then we have $\GG=\widetilde{\GG}$ because $\|\cdot\|_{\widetilde{\pi}}$ is equal to the original norm on $\C(\GG)$.
\item The quantum group $\widetilde{\GG}$ has, by construction, continuous counit. Moreover the comultiplication on $\C(\widetilde{\GG})$ obtained by extending that on $\Pol(\GG)$ is injective (cf.~the discussion in \cite{Dz}).
\item Suppose that $\GG\neq\widetilde{\GG}=\GGmax$. Note that it is obvious from the proof of Proposition \ref{mini} that the value of the counit (extended from $\Pol(\GG)$ to $\C(\widetilde{\GG})$) on the projection $p$ is $1$. Therefore by Propositions \ref{mini} and \ref{miprojT} we have that $\hGG$ has property $(\mathrm{T})$. 

We may therefore take for $\GG$ the reduced (minimal) version of a non-co-amenable compact quantum group whose dual does not have property $(\mathrm{T})$. Then 
\[
\GG\neq\widetilde{\GG}\neq\GGmax
\]
in the sense that the canonical morphisms $\C(\GGmax)\to\C(\widetilde{\GG})\to\C(\GG)$ are not isomorphisms.
\item The situation when $\GG\neq\widetilde{\GG}=\GGmax$ is also very interesting. We give and example of this phenomenon in Section \ref{exo}.
\end{enumerate}
\end{remark}

\section{Exotic quantum group norms}\label{exo}

As before we consider a discrete quantum group $\hGG$. In this section we will assume that $\hGG$ is infinite ($\dim\c0(\hGG)=\infty$) and that $\hGG$ has property $(\mathrm{T})$. Throughout this section we let $\Pi$ be the representation of $\C(\GGmax)$ defined as the direct sum of all its infinite dimensional irreducible representations. The corresponding corepresentation of $\hGG$ will be denoted by $\VV$:
\[
\VV=(\id\tens\Pi)\ww.
\]

\begin{proposition}
We have $\lambda\preccurlyeq\Pi$.
\end{proposition}

\begin{proof}
By Corollary \ref{notweakly} $\lambda$ does not weakly contain any finite dimensional representation. This means that the support of $\lambda$ in $\Spec\bigl(\C(\GGmax)\bigr)$ does not have any finite dimensional representations in its closure. Therefore the support of $\lambda$ is contained in the set of infinite dimensional representations, i.e.~the support of $\Pi$. Hence $\lambda\preccurlyeq\Pi$.
\end{proof}

The above result yields the following corollary:

\begin{corollary}
The seminorm $\|\cdot\|_\Pi$ defined on $\Pol(\GG)$ by $\Pi$ is a norm and $\C(\GGmin)$ is a quotient of the completion $\C(\GG_\Pi)$ of $\Pol(\GG)$ in the norm $\|\cdot\|_\Pi$.
\end{corollary}

\begin{theorem}\label{cool}
$\|\cdot\|_\Pi$ is a quantum group norm on $\Pol(\GG)$.
\end{theorem}

\begin{proof}
We will show that $(\Pi\tens\Pi)\comp\Delta\preccurlyeq\Pi$. We will do this using the language of corepresentations of $\hGG$ instead of that of representations of $\C(\GGmax)$.

Clearly it is enough to show that $\VV\tp\VV$ does not weakly contain a finite dimensional corepresentation. Assume the contrary and let $U$ be a finite dimensional corepresentation of $\hGG$ such that $U\preccurlyeq\VV\tp\VV$. By Remark \ref{fin} (and the fact that finite dimensional corepresentations decompose into direct sums of irreducible ones) we have $U\leq\VV\tp\VV$. Therefore, by Lemma \ref{little}
\[
U\tp{U^\cc}\leq(\VV\tp\VV)\tp(\VV\tp\VV)^\cc.
\]
But $U\tp{U^\cc}$ contains the trivial corepresentation of $\hGG$ (Remark \ref{bG}), so
\[
(\VV\tp\VV)\tp(\VV\tp\VV)^\cc\approx\VV\tp(\VV\tp\VV^\cc\tp\VV^\cc)
\]
(cf.~\cite[Formula (3.7)]{mu2}) contains the trivial corepresentation. Since $\hGG$ is unimodular (Theorem \ref{ft}\eqref{ftu}) we can use Theorem \ref{T}\eqref{UVW2} to see that $\VV$ must then contain a finite dimensional corepresentation. This is a contradiction with the construction of $\VV$.
\end{proof}

For our infinite, discrete property $(\mathrm{T})$ quantum group $\hGG$ we now obtain the following two results

\begin{corollary}
The compact quantum group $\GG_\Pi$ obtained via completion of $\Pol(\GG)$ in the norm $\|\cdot\|_\Pi$ does not admit a continuous counit.
\end{corollary}

This is evident, because our assumptions on $\hGG$ imply that $\eps\not\preccurlyeq\Pi$.

The same technique as the one used in the proof of Theorem \ref{cool} gives an answer to a question asked in \cite[End of Section 3]{BMT01}, namely if every \cst-norm on $\Pol(\GG)$ defined by a representation which weakly contains the regular one is a quantum group norm.

\begin{corollary}
Assume that $\hGG$ has a non-trivial finite-dimensional corepresentation $U$. Let $U_0$ be an irreducible subrepresentation of $U$ and let $\pi_0$ be the corresponding representation of $\C(\GGmax)$:
\[
U_0=(\id\tens\pi_0)\ww.
\]
Then the representation defined as the direct sum of all irreducible representations of $\C(\GGmax)$ except the trivial one weakly contains the regular representation and the associated norm is not a quantum group norm.
\end{corollary}

Let us now discuss one special case when the compact quantum group $\GG_\Pi$ has quite unexpected properties. Let us consider the cocommutative example with $\hGG=\Gamma$, a discrete Kazhdan group which is \emph{minimally almost periodic,} i.e.~it has no non-trivial finite dimensional irreducible representations.\footnote
{
Examples of discrete property $(\mathrm{T})$ groups which are minimally almost periodic have been constructed by Gromov in \cite{gr} (cf.~\cite[Theorem 3.4]{qd}, more explicit examples have been constructed in \cite{cr}). The result of Gromov provides (uncountably many) pairwise non-isomorphic infinite discrete property $(\mathrm{T})$ torsion groups. In particular, they cannot contain a non-abelian free subgroup, so by Tits' alternative (\cite[Section 42]{dlH}) they cannot be \emph{linear,} i.e.~subgroups of $\mathrm{GL}(N,\mathbb{K})$ for a field $\mathbb{K}$ of characteristic $0$. Moreover, by \cite[Lemma 3.5]{qd}, these groups are simple. Since they are not linear, they must be minimally almost periodic.
}
Then it is easily seen that $\Pi\oplus\eps$ is weakly equivalent to the universal representation of $\Pol(\GG)=\CC[\Gamma]$. In other words,
\begin{equation}\label{GpiGmax}
\widetilde{\GG_\Pi}=\GGmax.
\end{equation}

\begin{remark}
\noindent
\begin{enumerate}
\item Let us note that if $\hGG$ is an infinite discrete property $(\mathrm{T})$ group with only one irreducible finite dimensional corepresentation, namely the trivial one, then the minimal projection $p\in\C(\GGmax)$ associated to this representation has a very peculiar property. Let $\Delta_\Pi$ be the comultiplication on $\C(\GG_\Pi)$ and let $\rho\colon\C(\GGmax)\to\C(\GG_\Pi)$ be the quotient map. Of course we have $\rho(p)=0$. Note further that $\rho$ is faithful on $\Pol(\GG)\subset\C(\GGmax)$ (e.g.~because the regular representation $\lambda$ factors through $\rho$).  If we denote by $\Deltamax$ the comultiplication on $\C(\GGmax)$ then we have
\[
(\rho\tens\rho)\comp\Deltamax=\Delta_\Pi\comp\rho
\]
so that $(\rho\tens\rho)\Deltamax(p)=0$. However, due to the decomposition
\[
\C(\GGmax)\cong\C(\GG_\Pi)\oplus\CC{p}
\]
we clearly have
\[
\ker{\rho\tens\rho}=(\CC{p}\tens{p})\oplus\bigl({p}\tens\C(\GG_\Pi)\bigr)\oplus\bigl(\C(\GG_\Pi)\tens{p}\bigr)
\subset\C(\GGmax)\atens\C(\GGmax). 
\]
This means that $\Deltamax(p)\in\C(\GGmax)\atens\C(\GGmax)$, but $p\not\in\Pol(\GG)$, as $\rho(p)=0$. This example provides an answer to the famous question whether any element of a \cst-algebra $\C(\GG)$ whose image under the coproduct is a finite sum of simple tensors must belong to $\Pol(\GG)$. The affirmative answer for compact quantum groups with faithful Haar measure was given by S.L.~Woronowicz in \cite[Theorem 2.6(2)]{cqg}. Our example shows that this is not the case if the dual of $\GG$ is minimally almost periodic with property $(\mathrm{T})$. A crucial fact here is that there actually exists a comultiplication on $\C(\GG_\Pi)$ or, in other words, that $\|\cdot\|_\Pi$ is a quantum group norm. Note also, that a similar argument applies if $\hGG$ has only \emph{finitely many} irreducible finite dimensional corepresentations.
\item The reader will have noticed that in fact the situation that $\widetilde{\GG_\Pi}=\GGmax$ is equivalent to $\C(\GGmax)$ having no irreducible finite dimensional representations except $\eps$. In other words \eqref{GpiGmax} holds if and only if $\hGG$ is minimally almost periodic.
\end{enumerate}
\end{remark}

The above example leads to an important question, namely whether we have $\GGmin=\GG_\Pi$. It seems that this could actually be the case in some examples, but we have not been able to produce one (nor find it in literature). However, as the next proposition says, at least for the cocommutative examples, the case that $\GGmin\neq\GG_\Pi$ is rather common.

\begin{proposition}\label{Valette}
Let $\Gamma$ be an infinite discrete group with Kazhdan's property $(\mathrm{T})$ such that the regular representation of $\Gamma$ is weakly equivalent to the sum of all infinite dimensional irreducible representations of $\Gamma$. Then any non-amenable subgroup of $\Gamma$ must have finite index. In particular, $\Gamma$ cannot be linear.
\end{proposition}

Before proving this proposition let us state one lemma.

\begin{lemma}\label{Valette2}
Let $L$ be a subgroup of a discrete countable group $G$ such that the permutation representation $\lambda_{G/L}$ of $G$ on $\ell^2(G/L)$ is weakly contained the regular representation $\lambda_G$ of $G$. Then $L$ is amenable.
\end{lemma}

\begin{proof}
The characteristic function of $L$ is a positive definite function associated with the permutation representation $\lambda_{G/L}$ (consider the coefficient of $\lambda_{G/L}$ arising from the vector in $\ell^2(G/L)$ which is the
delta-function in the point $L$ of $G/L$). Therefore, since $\lambda_{G/L}$ is weakly contained in $\lambda_G$, the characteristic function of $L$ is a pointwise limit of positive definite functions with finite support. Restricting these functions to $L$ yields a net of finitely supported positive definite functions on $L$ approximating pointwise the constant function $1$. This proves that $L$ is amenable.
\end{proof}

\begin{proof}[Proof of Proposition \ref{Valette}]
Let $\Lambda$ be a non-amenable subgroup of $\Gamma$. Then, by Lemma \ref{Valette2}, the permutation representation $\lambda_{\Gamma/\Lambda}$ cannot be weakly contained in $\lambda_\Gamma$. By assumptions on $\Gamma$ there must be a finite dimensional representation $\sigma$ of $\Gamma$ weakly contained in $\lambda_{\Gamma/\Lambda}$ (there must be
an irreducible representation $\sigma$ weakly contained in $\lambda_{\Gamma/\Lambda}$ and not in $\lambda_\Gamma$, but all infinite dimensional ones are weakly contained in $\lambda_\Gamma$). Consider now the tensor product of $\sigma$ and its contragredient representation. This is weakly contained in the tensor product of
$\lambda_{\Gamma/\Lambda}$ with its contragredient which is equivalent to the tensor square of $\lambda_{\Gamma/\Lambda}$. Of course $\sigma\tens\sigma^\cc$ contains the trivial representation, so the square of $\lambda_{\Gamma/\Lambda}$ contains (strongly - by property $(\mathrm{T})$) the trivial representation. Now the tensor square of $\lambda_{\Gamma/\Lambda}$ is equivalent to the permutation representation of $\Gamma$ on $\ell^2((\Gamma/\Lambda)\times(\Gamma/\Lambda))$ with diagonal action.
If this representation has a fixed vector, then $\Gamma$ must have a finite orbit for the diagonal action on $(\Gamma/\Lambda)\times(\Gamma/\Lambda)$. (If a group $\Gamma$ acts on a set $S$ and the associated permutation representation in $\ell^2(S)$ has a non-zero fixed vector $\xi$, then expanding this vector in the canonical orthonormal basis and acting on it shows that the coefficients of $\xi$ are constant along orbits - therefore there
must be a finite orbit.) This means that $\Lambda$ has finite index in $\Gamma$ because if $(\gamma\Lambda,\gamma'\Lambda)$ is an element of $(\Gamma/\Lambda)\times(\Gamma/\Lambda)$ which has finite orbit, then there are $\gamma_1,\dotsc,\gamma_n,\gamma'_1,\dotsc,\gamma'_n\in\Gamma$ such that
\[
(x\gamma\Lambda,x\gamma'\Lambda)\in\bigl\{
(\gamma_1\Lambda,\gamma'_1\Lambda),\dotsc,(\gamma_n\Lambda,\gamma'_n\Lambda)\bigr\}
\]
for all $x\in\Gamma$. But $\{x\gamma\Lambda\st{x}\in\Gamma\}$ is all of $\Gamma/\Lambda$, so $\Gamma/\Lambda$ is contained in the union
\[
\gamma_1\Lambda\cup\dotsm\cup{\gamma_n\Lambda}.
\]
This establishes that any non-amenable subgroup of $\Gamma$ has finite index.

Since $\Gamma$ has property $(\mathrm{T})$ it is finitely generated, so if $\Gamma$ furthermore were linear the Tits alternative (\cite[Section 42]{dlH}) implies that it is either virtually solvable (which is impossible because it is non-amenable) or contains a non-abelian free subgroup. It is easy to see that then $\Gamma$ must also contain non-amenable subgroups of infinite index.
\end{proof}

It follows from Proposition \ref{Valette} that if we take $\hGG=\Gamma$ to be a linear infinite Kazhdan group, say $\Gamma=\mathrm{SL}(3,\mathbb{Z})$, admitting non-trivial finite-dimensional irreducible representations, then we have
\[
\GGmin\neq\GG_\Pi\neq\widetilde{\GG_\Pi}\neq\GGmax.
\]
\subsection*{Acknowledgments}
The authors are grateful to Vadim Alekseev, Narutaka Ozawa, Andrzej \.Zuk and, in particular, Alain Valette for helpful comments regarding the results in Section 9. The proof of Proposition \ref{Valette} is due to Alain Valette.

\end{document}